\documentclass[10pt,draft]{amsart}
\usepackage{amsmath,amssymb,amsthm}
\usepackage{amsmath,amssymb,amsthm}
\begin{document}
\newtheorem{lem}{Lemma}[section]
\newtheorem{prop}{Proposition}[section]
\newtheorem{cor}{Corollary}[section]
\numberwithin{equation}{section}
\newtheorem{thm}{Theorem}[section]

\theoremstyle{remark}
\newtheorem{example}{Example}[section]
\newtheorem*{ack}{Acknowledgments}

\theoremstyle{definition}
\newtheorem{definition}{Definition}[section]

\theoremstyle{remark}
\newtheorem*{notation}{Notation}
\theoremstyle{remark}
\newtheorem{remark}{Remark}[section]

\newenvironment{Abstract}
{\begin{center}\textbf{\footnotesize{Abstract}}%
\end{center} \begin{quote}\begin{footnotesize}}
{\end{footnotesize}\end{quote}\bigskip}
\newenvironment{nome}

{\begin{center}\textbf{{}}%
\end{center} \begin{quote}\end{quote}\bigskip}

\newcommand{\triple}[1]{{|\!|\!|#1|\!|\!|}}

\newcommand{\xx}{\langle x\rangle}
\newcommand{\ep}{\varepsilon}
\newcommand{\al}{\alpha}
\newcommand{\be}{\beta}
\newcommand{\de}{\partial}
\newcommand{\la}{\lambda}
\newcommand{\La}{\Lambda}
\newcommand{\ga}{\gamma}
\newcommand{\del}{\delta}
\newcommand{\Del}{\Delta}
\newcommand{\sig}{\sigma}
\newcommand{\ome}{\omega}
\newcommand{\Ome}{\Omega}
\newcommand{\C}{{\mathbb C}}
\newcommand{\N}{{\mathbb N}}
\newcommand{\Z}{{\mathbb Z}}
\newcommand{\R}{{\mathbb R}}
\newcommand{\Rn}{{\mathbb R}^{n}}
\newcommand{\Rnu}{{\mathbb R}^{n+1}_{+}}
\newcommand{\Cn}{{\mathbb C}^{n}}
\newcommand{\spt}{\,\mathrm{supp}\,}
\newcommand{\Lin}{\mathcal{L}}
\newcommand{\SSS}{\mathcal{S}}
\newcommand{\F}{\mathcal{F}}
\newcommand{\xxi}{\langle\xi\rangle}
\newcommand{\eei}{\langle\eta\rangle}
\newcommand{\xei}{\langle\xi-\eta\rangle}
\newcommand{\yy}{\langle y\rangle}
\newcommand{\dint}{\int\!\!\int}
\newcommand{\hatp}{\widehat\psi}
\renewcommand{\Re}{\;\mathrm{Re}\;}
\renewcommand{\Im}{\;\mathrm{Im}\;}

\title{Existence of maximizers for Sobolev-Strichartz inequalities}

\author{Luca Fanelli}
\address{Luca Fanelli: Universidad del Pais Vasco, Departamento de
Matem$\acute{\text{a}}$ticas, Apartado 644, 48080, Bilbao, Spain}
\email{luca.fanelli@ehu.es}

\author{Luis Vega}
\address{Luis Vega: Universidad del Pais Vasco, Departamento de
Matem$\acute{\text{a}}$ticas, Apartado 644, 48080, Bilbao, Spain}
\email{luis.vega@ehu.es}

\author{Nicola Visciglia}
\address{Nicola Visciglia: Universit\`a di Pisa, Dipartimento di Matematica, Largo B.
Pontecorvo 5, 56100 Pisa, Italy}
\email{viscigli@dm.unipi.it}

\subjclass[2000]{35J10, 35L05.}

\keywords{Fourier restriction theorems, Strichartz estimates}

\maketitle

\begin{abstract}
We prove the existence of maximizers
of Sobolev-Strichartz estimates for a general class of propagators,
involving relevant examples, as for instance the wave, Dirac and the hyperbolic Schr\"odinger flows.
\end{abstract}

\section{Introduction}\label{sec:intro}
In the recent years, an increasing interest has 
been devoted to the problem of existence of maximizers for Strichartz inequalities, 
and more in general Fourier Restriction Theorems. \\
Recall that, given two Banach spaces
$(X, |.|_X)$, $(Y, |.|_Y)$,
and a linear and continuos operator $T\in {\mathcal L}(X, Y)$,
it is customary to define maximizer
for $T$ any $x_0\in X$ such that
\begin{equation}
\|x_0\|_X=1,
\qquad
 \|Tx_0\|_Y=\sup_{\|x\|_X=1} \|Tx\|_Y.
\end{equation}
Notice that the existence (as well as the definiton) 
of maximizerz 
depends on the specific norm introduced on $X$ and $Y$.\\
In this paper we consider the existence of maximizers
for a large class of propagators of the form $e^{ith(D)}$, 
and we focus our attention to cases in which some 
Sobolev-Strichartz estimates of the following type hold
\begin{equation}\label{eq:0}
  \left\|e^{ith(D)}f\right\|_{L^r_{t,x}}\leq C\|f\|_{\dot H^s},
  \qquad
  s>0,\ \ r>2.
\end{equation}
where in general $f$ is a vector valued function $f=(f_1,...,f_n):\R^d\rightarrow \C^n$
and
\begin{equation}\label{eq:Hsobolev}
\|f\|_{\dot H^s}^2=\sum_{j=1}^n \|f_j\|_{\dot H^s}^2.\end{equation}
Concerning the definition of the norm $L^{r}_{t,x}$ in the l.h.s. of \eqref{eq:0}
recall that in general for any given norm $|.|_{\C^n}$ on $\C^n$
one can define the corresponding mixed Lebesgue norms  $L^p_tL^q_x$  
for vector valued functions $F=F(t,x):\R^{1+d}\to\C^n$ as follows
\begin{equation*}
  \| F(t,x)\|_{L^p_tL^q_x}^p= \int \left (\int |F(t,x)|^q_{\C^n} dx\right )^\frac pq dt,
\end{equation*} 
where $1\leq p, q\leq \infty$. In the case $p=q=r$ we shall write shortly
$$L^r_tL^r_x=L^r_{t,x}.$$
%In the sequel, we also treat the case of systems, 
%in which we generally deal with vector-valued functions $F=F(t,x):\R^{1+d}\to\C^n$. 
%Notice that each norm on $\C^n$ induces in a natural way 
%$L^p_tL^q_x$-norm on functions with values in 
%$\C^n$. 
\begin{remark}
Recall that all the norms on $\C^n$ are equivalent
and hence also the corresponding norms
induced on the spaces $L^p_tL^q_x$ are equivalent. Due to this fact 
in \eqref{eq:0} it is not necessary to specify
the norm on $\C^n$ with respect to which we are working with
(of course provided that the corresponding constant $C$
is suitably modified).
On the other hand the existence of maximizers for the inequality \eqref{eq:0}
could be 
affected in principle by changing the corresponding norm on $\C^n$.
\end{remark}
Estimates of the kind \eqref{eq:0}  with a loss of derivatives with respect 
to the datum $f$, are natural in some relevant cases, as for instance the wave and Dirac equations. 
Let us point our attention to the Cauchy problem for the following system
\begin{equation}\label{eq:eq}
  \begin{cases}
    i\partial_t u+h(D) u=0
    \\
    u(0,x)=f(x),
  \end{cases}  
\end{equation}
where $u(t,x)=(u_1(t,x),\dots,u_n(t,x)):\R^{1+d}\to\C^n$, $f(x)=(f_1(x),\dots,f_n(x)):\R^d\to\C^n$. 
Here we denote by $h(D)=\mathcal F^{-1}\left(h(\xi)\mathcal F\right)$, where 
$\mathcal F$ is the standard Fourier transform; in addition, the symbol $h(\xi)\in\mathcal M_{n
\times n}(\C)$, $n\geq1$, 
is assumed to be a matrix-valued function $h(\xi)=\left(h_{ij}(\xi)\right)_{ij=1,\dots n}$. 
In the sequel, we always make the following abstract assumptions:
\begin{itemize}
  \item[(H1)]
  there exists $0<s<\frac d2$ such that 
  the problem \eqref{eq:eq} is globally well-posed in $\dot H^s$, and 
the unique solution is given via propagator $u(t,x)=e^{ith(D)}f(x)$;
  \item[(H2)]  the flow $e^{ith(D)}$ is unitary onto $\dot H^s$, i.e.
  \begin{equation*}
    \left\|e^{ith(D)}f\right\|_{\dot H^s}=\|f\|_{\dot H^s}
    \qquad
    \forall t\in\R,
  \end{equation*}
  where $s$ is the same as in (H1), and $\|\cdot\|_{\dot H^s}$ is defined in \eqref{eq:Hsobolev}.
\end{itemize}

%We always denote by
%\begin{equation*}
%  \|f\|_{\dot H^s}^2=\sum_{i=1}^n \|f_i\|_{\dot H^s}^2.
%\end{equation*}

%The three relevant examples we have in mind are the following:
%\begin{itemize}
%  \item
%  the wave propagator $h(D)=|D|$ with $N=1$;
%  \item
%  the Dirac propagator: $h(D)=-i\sum_{j=1}^3\alpha_j\partial_j$, where $\alpha_1,\alpha_2,\alpha_3
%  \in\mathcal M_{4\times4}(\C)$ are the Dirac matrices in the standard Pauli representation. Here we have $N=4$;
%  \item
%  the hyperbolic Schr\"odinger propagator $h(D)=\sum_{j=1}^r\partial_j^2-\sum_{j=r+1}^d
%  \partial_j^2$, where $d\geq2$ and $1\leq r<d$. In this case we have $N=1$.
%\end{itemize}

In order to introduce our problem,
assume for the moment that $h_{i,j}(\xi)$ is homogenous of some degree $k>0$
for every $i, j\in \{1,..,n\}$, namely
\begin{equation}\label{eq:omog}
  h_{ij}(\xi):\R^d\to\C,
  \qquad
  h(\lambda\xi)=\lambda^kh(\xi)
  \qquad
  \forall\lambda>0.
\end{equation} 
Then, equation \eqref{eq:eq} is invariant under the scaling 
$v_\lambda(t,x)=u(\lambda^kt,\lambda x)$; as a consequence, if a Strichartz estimate of the following type holds
\begin{equation}\label{eq:stri}
\left\|e^{ith(D)} f\right\|_{L^p_tL^q_x}\leq C\| f\|_{\dot H^{s}},
\end{equation}
for some $s>0$ and some constant $C>0$, then $p$ and $q$ have to satisfy the following scaling condition
\begin{equation}\label{eq:admis}
\frac kp+\frac dq=\frac d2-s.
\end{equation}
Moreover, by the Sobolev embedding $\dot H^s\subset L^{\frac{2d}{d-2s}}$, 
for $0<s<d/2$, and the $\dot H^s$-preservation
\begin{equation*}
  \left\|e^{ith(D)}f\right\|_{L^\infty_t\dot H^s_x}=\|f\|_{\dot H^s}
\end{equation*}
(assumption (H2) above),
we get
\begin{equation}\label{eq:1}
  \left\|e^{ith(D)}f\right\|_{L^\infty_tL^{\frac{2d}{d-2s}}_x}\leq C\|f\|_{\dot H^s},
  \qquad
  0<s<\frac d2,
\end{equation}
for some constant $C>0$.
Hence, if an estimate of the type \eqref{eq:stri} holds with $p<q$, then by interpolation with \eqref{eq:1} we obtain
\begin{equation}\label{eq:est}
  \left\|e^{ith(D)}f\right\|_{L^r_{t,x}}\leq C\|f\|_{\dot H^s},
  \qquad
  0<s<\frac d2,
  \qquad
  r=\frac{2(k+d)}{d-2s}.
\end{equation}
The aim of this paper is to prove that, as soon as an estimate 
as \eqref{eq:est} holds, with a strictly positive $s>0$, 
then the best constant of the estimate is achieved by some maximizing functions, 
independently on the norm which is fixed on the target 
$\C^n$ (and consequently on the corresponding definition of the $L^r_{tx}$-norm). \\
This kind of problem has been recently studied by several authors, who treated 
separtely different propagators. We firstly mention Kunze \cite{K}, who proved the existence 
of maximizers of the $L^6_{t,x}$-Strichartz inequality for the 1D Schr\"odinger propagator. 
Later, Foschi \cite{F} succeeded in characterizing the best constant of the inequality 
and also the shape of the maximizers, for the 1D and 2D-Schr\"odinger propagators. 
In the same paper, the author treats the $L^6_{t,x}$ (in 2D) and the $L^4_{t,x}$ (in 3D) 
Strichartz estimates for the wave equation, which hold to the scale of 
$\dot H^{\frac12}\times\dot H^{-\frac12}$-initial data. Recently, Bez and Rogers \cite{BR} 
computed the best constant and described the shape of maximizers of the 
$L^4_{t,x}$-Strichartz estimate for the wave equation, with initial 
data in the energy space $\dot H^1\times L^2$, in space dimension $d=5$. 
We also mention the papers by Shao \cite{S} for the Sch\"odinger
equation and by Bulut \cite{B} for the wave equation, in which the existence of 
maximizers for anisotropic Strichartz 
estimates in spaces 
$L^p_tL^q_x$, for general couples $(p,q)$, is also proved. Finally, Ramos 
\cite{R} proves that there exist maximizers for the isotropic Strichartz estimates for the wave equation, 
at the scale $\dot H^{\frac12}\times\dot H^{-\frac12}$.

Here we give a unified (and simple) proof of the existence of maximizers of 
\eqref{eq:est}, when $s>0$, which involves a large class of examples of propagators.
Our main result is the following.
\begin{thm}\label{thm:main}
  Let assumptions (H1), (H2) be satisfied for some $0<s<\frac d2$, and let $h(\xi)$  
satisfy \eqref{eq:omog} for some $k>0$. Moreover, assume that, for some $2\leq p<q\leq\infty$  
  \begin{equation}\label{eq:i}
    \left\|e^{ith(D)} 
f\right\|_{L^p_tL^q_x}\leq C\| f\|_{\dot H^s},
  \end{equation}
  so that, for $r=\frac{2(k+d)}{d-2s}$, we also have 
  \begin{equation}\label{eq:estr}
    \left\|e^{ith(D)} f\right\|_{L^r_{t,x}}\leq M\| f\|_{\dot H^s},
  \end{equation}
  with
  \begin{equation}\label{eq:M}
    M:=\sup_{\|f\|_{\dot H^s}=1}\left\|e^{ith(D)} f\right\|_{L^r_{t,x}}.
  \end{equation}
  Then there exists $f_0\in\dot H^s$ such that 
  \begin{equation}\label{eq:tesi}
  \|f_0\|_{\dot H^s}=1,
  \qquad
    \left\|e^{ith(D)}f_0\right\|_{L^r_{t,x}}=M.
  \end{equation}
\end{thm} 
\begin{remark}
We underline that in Theorem \ref{thm:main} no specific norm  is fixed 
on the target $\C^n$ (and hence also on the corresponding $L^r_{t,x}$ norm). Neverthless
we prove in general the existence of maximizers for the Strichartz estimate \eqref{eq:estr}. 
We also underline that 
in Theorem \ref{thm:main} we assume that the norm on $\dot H^s$
is the one defined in \eqref{eq:Hsobolev}. Indeed the Hilbert structure
of the norm $\|.\|_{\dot H^s}$ is crucial in our argument.
\end{remark}
\begin{remark}\label{rem:r2}
  Notice that the conditions $2\leq p<q\leq\infty$ imply that 
$r>2$, which will be crucial in the proof of Theorem \ref{thm:main}.
\end{remark}
\begin{remark}\label{rem:reg}
  As we see in the sequel, the proof of Theorem \ref{thm:main} is quite simple, and it is based on a suitable 
variant of a classic result by Br\'ezis and Lieb (\cite{BL}, \cite{L}, see Section \ref{sec:proof} below). 
On the other hand, the technique we use does not allow 
us to consider the case $s=0$ for which, in fact, some extra ingredients 
(being typically suitable improvements of Strichartz estimates) 
are needed, as substitutes of \eqref{eq:i}.
  \end{remark}
We shall now give some examples of applications of the previous Theorem.
\begin{example}[Wave propagator]\label{ex:wave}
The Strichartz estimates for the wave propagator $e^{it|D|}$ 
(\cite{GV}, \cite{KT}), in dimension $d\geq2$, are the following:
\begin{equation}\label{eq:striwave}
  \left\|e^{it|D|}f\right\|_{L^p_tL^q_x}\leq C\|f\|_{\dot H^{\frac 1p-\frac 1q+\frac 12}},
\end{equation}
under the admissibility condition 
\begin{equation}\label{eq:admiswave}
  \frac 2p+\frac{d-1}q=\frac{d-1}2,
  \qquad
  p\geq2.
  \qquad
  (p,q)\neq(2,\infty).
\end{equation}
In this case, the gap of derivatives $\frac 1p-\frac 1q+\frac 12\geq0$ is 
null only in the case of the energy estimate $(p,q)=(\infty,2)$. 
In particular, we have 
\begin{equation}\label{eq:striwave2}
  \left\|e^{it|D|}f\right\|_{L^{\frac{2(d+1)}{d-1}}_{t,x}}\leq C\|f\|_{\dot H^{\frac12}}
  \qquad
  d\geq2,
\end{equation}
which is in fact the original estimate proved by Strichartz in \cite{STR}.
More generally, by Sobolev embedding one also obtains that
\begin{equation}\label{eq:striwave20}
  \left\|e^{it|D|}f\right\|_{L^{\frac{2(d+1)}{d-1-2\sigma}}_{t,x}}\leq C\|f\|_{\dot H^{\frac12+\sigma}},
  \qquad
  0\leq\sigma
<\frac{d-1}2,
\qquad
  d\geq2.
\end{equation}
Theorem \ref{thm:main} gives a short and simple proof of the fact that, for any dimension $d\geq2$, 
and $0<\sigma<\frac{d-1}2$, the best constant in 
\eqref{eq:striwave20} is achieved on some maximizing function $f$.
\end{example}
\begin{example}[Dirac equation]\label{ex:dirac}
  Consider the (massless) Dirac operator 
  \begin{equation*}
   \mathcal D:=\frac1i\sum_{j=1}^3\alpha_j\partial_j,
  \end{equation*} 
   which is defined in dimension $d=3$. 
  Here $\alpha_1,\alpha_2,\alpha_3\in\mathcal M_{4\time4}(\C)$ are the so called {\it Dirac matrices}, 
which are $4\times4$-hermitian matrices, $\alpha_j^t=\overline{\alpha_j}$, $j=1,2,3$. 
They are defined as
  \begin{equation*}
    \alpha_1=
    \left(
    \begin{array}{cccc}
    0 & 0 & 0 & 1
    \\
    0 & 0 & 1 & 0
    \\
    0 & 1 & 0 & 0
    \\
    1 & 0 & 0 & 0
    \end{array}
    \right),
    \ \ \ 
    \alpha_2=
    \left(
    \begin{array}{cccc}
    0 & 0 & 0 & -i
    \\
    0 & 0 & i & 0
    \\
    0 & -i & 0 & 0
    \\
    i & 0 & 0 & 0
    \end{array}
    \right),
    \ \ \ 
    \alpha_3=
    \left(
    \begin{array}{cccc}
    0 & 0 & 1 & 0
    \\
    0 & 0 & 0 & -1
    \\
    1 & 0 & 0 & 0
    \\
    0 & -1 & 0 & 0
    \end{array}
    \right)
  \end{equation*}
  or equivalently $\alpha_j=\left(\begin{array}{cc}0 & \sigma_j \\ \sigma_j & 0\end{array}\right)$, where $\sigma_j$ is the $j^{\text{th}}$ $2\times2$-Pauli matrix, $j=1,2,3$.
  
  Since $\mathcal D^2=-\Delta I_{4\times4}$, the Strichartz estimates for the massless Dirac operator are the same as for the 3D wave equation (see \cite{DF}):  \begin{equation}\label{eq:stridir}
    \left\|e^{it\mathcal D} f\right\|_{L^p_tL^q_x}\leq C\|
f\|_{\dot H^{\frac 1p-\frac 1q+\frac 12}},
  \end{equation}
  with the admissibility condition
  \begin{equation}\label{eq:admisdir}
  \frac 2p+\frac{2}q=1,
  \qquad
  p>2.
\end{equation}
In particular we have
\begin{equation}\label{eq:stridir2}
  \left\|e^{it\mathcal D} f\right\|_{L^{\frac{8}{2-2\sigma}}_{t,x}}\leq C\| f\|_{\dot H^{\frac12+\sigma}},
  \qquad
  0\leq\sigma<1.
\end{equation}
Hence also in this case a loss of derivatives with respect to the initial datum is natural in the estimate.
Also in this case Theorem \ref{thm:main} applies, and it proves that there exist maximizers for \eqref{eq:stridir2}, in the range $0<\sigma<1$, which at our knowledge is not a known fact.
\end{example}
\begin{example}[Hyperbolic Schr\"odinger equation]\label{ex:schro}
  Let us now consider the hyperbolic Schr\"odinger operator 
  $L:=\sum_{j=1}^m\partial_j^2-\sum_{j=m+1}^d
  \partial_j^2$, for $d\geq2$ and $1\leq m<d$. In this case, Strichartz estimates are the same as for the Schr\"odinger propagator, namely
  \begin{equation}\label{eq:strischro}
     \left\|e^{itL}f\right\|_{L^p_tL^q_x}\leq C\|f\|_{\dot H^s},
  \end{equation}
  with the admissibility condition
  \begin{equation}\label{eq:admisschro}
    \frac 2p+\frac dq=\frac d2-s,
    \qquad
    p\geq2,
    \qquad
    (p,q)\neq(2,\infty).
  \end{equation}
  In particular, one has
  \begin{equation}\label{eq:strischro2}
    \left\|e^{itL}f\right\|_{L^{\frac{2(d+2)}{d-2s}}_{t,x}}\leq C\|f\|_{\dot H^s},
  \end{equation}
  for any $0\leq s<\frac d2$. The only case in which the existence of maximizers for 
  \eqref{eq:strischro2} is known is $d=2,s=0$. Indeed, the result by Rogers and Vargas in \cite{RV} contains all the ingredients which are necessary to prove the profile decomposition for bounded sequences in $L^2$, with respect to the propagator $e^{itL}$, in dimension $d=2$; the existence of maximizers follows from this fact by a general argument in the spirit of \cite{S}.
  It is a matter of fact that Theorem \ref{thm:main} applies for any $d\geq2$ and $0<s<\frac d2$, but we remark that it cannot include the case $s=0$. We finally remark, our argument is rather simple and does not involve any profile decomposition-theorems associated to the propagator.
\end{example}
  In the statement of Theorem \ref{thm:main}, the homogeneity assumption \eqref{eq:omog} is required; in fact, this is just put in order to write the explicit dependence of $r$ on $k,s$ in \eqref{eq:estr}. Motivated by the case of the wave equation, in which the homogeneity property does not hold (see example \ref{rem:wave} below), we now state a more general version of the previous theorem.

Let us first introduce the notations
\begin{equation*}
  \dot H^s=\dot H^{s_1}\times\dots\times \dot H^{s_n}
  \qquad
  s:=(s_1,\dots,s_n),
\end{equation*}  
\begin{equation}\label{eq:sobole}
  \|f\|_{\dot H^s}^2=\sum_{i=1}^n \|f_i\|_{\dot H^{s_i}}^2
\end{equation}
where $f=(f_1,\dots,f_n)\in \dot H^s$.
Consider equation \eqref{eq:eq}; in addition to (H1) and (H2)
(where the norm $\dot H^s$ in this case is the more general 
one defined in \eqref{eq:sobole})  assume:
\begin{itemize}
\item[(H3)] 
there exists $j\in \{1,..,n\}$ and $2\leq p<q\leq \infty$ such that, 
if $(u_1,...,u_n)$ solves \eqref{eq:eq} then 
\begin{align*}
\|u_j\|_{L^{p}_{t}L^q_x} & \leq C_j \|f \|_{\dot H^s};
\end{align*}
\item[(H4)]  
if $u=(u_1,,, u_n)$ solves \eqref{eq:eq}, then 
\begin{equation*}
  u_\lambda= \left(\lambda^{\frac d2 -s_1} u_1(\lambda t, \lambda x),\dots,
\lambda^{\frac d2-s_n} u_n(\lambda t, \lambda x)\right) 
\end{equation*}
solves \eqref{eq:eq}, 
for any $\lambda>0$. 
\end{itemize}
Notice that, by a scaling argument, we have
\begin{equation*}
   \frac 1p + \frac dq= \frac d2 - s_j.
\end{equation*}
Moreover, interpolating with the energy estimate one also obtains 
\begin{equation}\label{eq:estrappendix}
    \left\|u_j \right\|_{L^r_{t,x}}\leq M\| f\|_{\dot H^s},
\end{equation}
 where
  $r=\frac{2(1+d)}{d-2s_j}$.
We have the following result.
\begin{thm}\label{thm:2}
Assume the operator $h(D)$ is such that (H1), (H2), (H3), (H4) 
are satisfied for some  $s_1,\dots,s_n\in \R$, $2\leq p<q\leq \infty$, and $0<s_j<\frac d2$, where 
$j$ is the one in (H3). Consider estimate \eqref{eq:estrappendix} and let  
\begin{equation}\label{eq:Mappendix}
    M:=\sup_{\|f\|_{\dot H^s}=1}\left\| u_j\right\|_{L^r_{t,x}}.
  \end{equation}
 Then there exists $f_0\in\dot H^s$ such that 
  \begin{equation}\label{eq:tesi2}
  \|f_0\|_{\dot H^s}=1,
  \qquad
    \left\|v_j \right\|_{L^r_{t,x}}=M  
    \end{equation}
where $(v_1,\dots,v_n)=e^{ith(D)} f_0$
\end{thm}
\begin{remark}
  In other words, Theorem \ref{thm:2} states that the problem of existence of 
maximizers can be solved analogously for Sobolev-Strichartz estimates involving one single component of the solution of equation \eqref{eq:eq}. In a 
completely analogous way, one could treat 
the case of estimates for a selected group of components of the solution.
The proof of Theorem \ref{thm:2} is completely identical to the one of Theorem \ref{thm:main}: indeed it is sufficient to repeat the same argument for the composition operator $T_j(t)=\pi_j\circ e^{ith(D)}$, where $\pi_j:\C^N\to\C$ is the projection onto 
the $j$-th component. We will omit the details of the proof.
\end{remark}

We pass now to show the main application of Theorem \ref{thm:2}, on solutions of the wave equation.
\begin{example}[Wave equation]\label{rem:wave}
Consider now the wave equation
\begin{equation}\label{eq:WAVE}
\begin{cases}
  \partial_t^2u-\Delta u=0 \quad\text{in }\R^{1+d}
  \\
  u(0,x)=u_0(x)
  \quad
  \partial_tu(0,x)=u_1(x).
\end{cases}
\end{equation}
The vector-variable $V:=(u,\partial_tu)^t$ is uniquely given by $V(t,x)=e^{ith(D)}V_0(x)$, where $V_0=(u_0,u_1)^t$, and $h(D)=\left(
\begin{array}{cc}
0 & -i
\\
-i\Delta & 0
\end{array}
\right).$ 
The Strichartz estimates for \eqref{eq:WAVE} can be immediately deduced by \eqref{eq:striwave20}, writing the solution $u$ as 
\begin{equation*}
  u(t,\cdot)=\cos(t|D|)u_0(\cdot)+\frac{\sin(t|D|)}{|D|}u_1(\cdot);
\end{equation*}
in particular, in terms of the vector variable $V$, denoting by $f=(u_0,u_1)^t$, one has
\begin{equation}\label{eq:striwave200}
  \left\|u\right\|_{L^{\frac{2(d+1)}{d-1-2\sigma}}_{t,x}}\leq C\|f\|_{\dot H^{\frac12+\sigma}\times\dot H^{-\frac12+\sigma}},
  \qquad
  0\leq\sigma
<\frac{d-1}2,
\qquad
  d\geq2.
\end{equation}
Foschi (\cite{F}) proved that the best constant in \eqref{eq:striwave200} is achieved, in the cases 
$ d=2,3$, $\sigma=0$; moreover, he can characterize the shape of the maximizers. Later,
Bulut (\cite{B}) proved the same (also the anisotropic version of \eqref{eq:striwave200}, with $p\neq q$); she can treat the range $\frac12\leq \sigma<\frac{d-1}2$, with $d\geq3$. More recently, Bez and Rogers (\cite{BR}) can prove the existence of maximizers for \eqref{eq:striwave200}, and characterize their shape, in the case $d=5$, $\sigma=\frac12$. 

Now notice that Theorem \ref{thm:2} applies (with the choice $j=1$ in assumption (H3)).
It implies that there exist maximizers for \eqref{eq:striwave200}, for any dimension $d\geq2$ and any $0<\sigma<\frac{d-1}2$. We remark that Ramos (\cite{R}) proves a refinement of  
\eqref{eq:striwave200} in the case 
$\sigma=0$, which in particular implies the existence of maximizers, in this specific case of the wave equation. 

\end{example}

The rest of the paper is devoted to the proof of Theorem \ref{thm:main}.

\section{Proof of Theorem \ref{thm:main}}\label{sec:proof}
Before starting with the proof of our main theorem, we need to recall two fundamental results which will play a role in the sequel. The first one is a variant of a well known result obtained by Br\'ezis and Lieb in \cite{BL}, \cite{L}.
\begin{prop}[\cite{FVV}]\label{prop:CC}
Let $\mathcal H$ be a Hilbert space and
$T\in {\mathcal L}({\mathcal H}, L^p(\R^d))$
for a suitable $p\in (2,\infty)$.
Let $\{h_n\}_{n\in \N}\in \mathcal H$
such that:
\begin{align}
\|h_n\|_{\mathcal H} & = 1;
\label{eq:1000}
\\
\lim_{n\rightarrow \infty} \|T h_n\|_{L^p(\R^d)} & =\|T\|_{{\mathcal L}({\mathcal H}, L^p(\R^d))};
\label{eq:2000}
\\
h_n \rightharpoonup \bar h & \neq 0;
\label{eq:3000}
\\
T(h_n)\rightarrow T(\bar h) &\ \text{a.e. in } \R^d.
\label{eq:4000}
\end{align}
Then $h_n\rightarrow \bar h$ in $\mathcal H$,
in particular $\|\bar h\|_{\mathcal H}=1$ and
$\|T(\bar h)\|_{L^p(\R^d)}=\|T\|_{{\mathcal L}({\mathcal H}, L^p(\R^d))}$.
\end{prop}
\begin{remark}
The main difference between Proposition \ref{prop:CC}
and Lemma 2.7 in \cite{L} is
that we only need to assume
weak convergence in the Hilbert space $\mathcal H$ for
the maximizing sequence $h_n$.
On the other hand the argument in
\cite{L} works for operators defined between general
Lebesgue spaces and not necessarily
in the Hilbert spaces framework.
\end{remark}
Proposition \ref{prop:CC} has been proved in \cite{FVV}, in the scalar case. The proof in the vector-case is completely analogous and will be omitted. The second tool we need is a byproduct of a well known result by G\'erard in \cite{G}, in which the lack of compactness of the Sobolev embedding $\dot H^s\subset L^{\frac{2d}{d-2s}}$ is classified.
\begin{prop}[\cite{G}]\label{prop:gerard}
  Let $0<s<\frac d2$, and let $w_n\in\dot H^s$ be a sequence such that
  \begin{align}\label{eq:gerard1}
    \|w_n\|_{\dot H^s} & =1,
    \qquad
    n=1,2,\dots
    \\
    \inf_{n\geq1}\|w_n\|_{L_x^{\frac{2d}{d-2s}}} & 
    \geq \epsilon>0.
    \label{eq:gerard2}
  \end{align}
  Then there exist a sequence of parameters $\lambda_n>0$, a sequence of centers $x_n\in\R^d$, and a non-zero function $0\neq v\in\dot H^s$ such that
  \begin{equation}\label{eq:gerard3}
    v_n(x):=\lambda_n^{\frac d2-s}w_n(\lambda_n(x-x_n))\rightharpoonup v,
  \end{equation}
  weakly in $\dot H^s$, as $n\to\infty$.
\end{prop} 
\begin{remark}\label{rem:gerard}
  We remark that the condition $s>0$ is crucial in the previous proposition, which is in fact false in the case $s=0$. 
\end{remark}
We are now ready to perform the proof of our main theorem.
\subsection{Proof of Theorem \ref{thm:main}}
Let $M$ be as in \eqref{eq:M} and let $u_n\in\dot H^s$ be a maximizing sequence,  i.e.
\begin{equation}\label{eq:un}
  \|u_n\|_{\dot H^s}=1,
  \qquad
  \lim_{n\to\infty}\left\|e^{ith(D)}u_n\right\|_{L^r_{t,x}}=\left\|e^{ith(D)}\right\|_{\mathcal L\left(
  \dot H^s,L^r_{t,x}\right)}=M,
\end{equation}
with $r=\frac{2(k+d)}{d-2s}$. Our aim is to prove that, by a suitable remodulation of $u_n$, we can obtain a new maximizing sequence for which Proposition \ref{prop:CC} applies. In fact, since 
$u_n$ is uniformly bounded in $\dot H^s$, it admits a weak limit, which in principle could be zero.

Notice that, by Sobolev embedding,
\begin{equation}\label{eq:10}
  \left\|e^{ith(D)}u_n\right\|_{L^\infty_tL^{\frac{2d}{d-2s}}_x}
  \leq C
  \left\|e^{ith(D)}u_n\right\|_{L^\infty_t\dot H^s_x}=
  C\|u_n\|_{\dot H^s_x}=C,
\end{equation}
for some constant $C>0$. Moreover, by assumption \eqref{eq:i}, there exist 
$2\leq p<q\leq\infty$ such that
\begin{equation}\label{eq:20}
  \left\|e^{ith(D)}u_n\right\|_{L^p_tL^q_x}
  \leq C\|u_n\|_{\dot H^s_x}=C,
\end{equation}
for another constant $C>0$. Hence, by \eqref{eq:un}, \eqref{eq:20} and the H\"older inequality we can estimate,
for $n$ sufficiently large,
\begin{align*}
  \frac M2 
  \leq
  \left\|e^{ith(D)}u_n\right\|_{L^r_{t,x}}
  & 
  \leq
  \left\|e^{ith(D)}u_n\right\|_{L^\infty_tL^{\frac{2d}{d-2s}}_x}^{1-\frac pr}
  \left\|e^{ith(D)}u_n\right\|_{L^p_tL^q_x}^{\frac pr}
  \\
  &
  \leq
  C^{\frac pr}
  \left\|e^{ith(D)}u_n\right\|_{L^\infty_tL^{\frac{2d}{d-2s}}_x}^{1-\frac pr}
  \nonumber
\end{align*}
(where we have used \eqref{eq:admis}).
The last estimate implies that
\begin{equation}\label{eq:30}
  \left\|e^{ith(D)}u_n\right\|_{L^\infty_tL^{\frac{2d}{d-2s}}_x}
  \geq
  \left(\frac{M}{2C^{\frac pr}}\right)^{\frac{r}{r-p}}=:\epsilon>0,
\end{equation}
for $n$ sufficiently large. As a consequence, there exists a sequence of times $t_n\in\R$ such that
\begin{equation}\label{eq:40}
  \left\|e^{it_nh(D)}u_n\right\|_{L^{\frac{2d}{d-2s}}_x}\geq\frac\epsilon2>0,
\end{equation}
for any $n$ sufficienlty large. Now denote by
\begin{equation*}
  w_n(x)=e^{it_nh(D)}u_n(x),
\end{equation*}
and notice that $w_n$ is still a maximizing sequence, i.e. 
$\|e^{ith(D)} w_n\|_{L^r_{t,x}}\to M$, as $n\to\infty$. Moreover, we have
\begin{align}
\|w_n\|_{\dot H^s} & = \|u_n\|_{\dot H^s}=1
\label{eq:gerard11}
\\
\|w_n\|_{L^{\frac{2d}{d-2s}}_x} & \geq\frac\epsilon2>0;
\label{eq:gerard22}
\end{align}
hence, by Proposition \ref{prop:gerard}, there exist two sequences $\lambda_n>0$, $x_n\in\R^d$, and a non-zero function $0\neq v\in\dot H^s$, such that
\begin{equation}\label{eq:50}
  v_n(x):=\lambda_n^{\frac d2-s}w_n(\lambda_n(x-x_n))=
  \lambda_n^{\frac d2-s}e^{it_nh(D)}u_n(\lambda_n(x-x_n))\rightharpoonup v\neq0
\end{equation}
weakly in $\dot H^s$, as $n\to\infty$. By scaling, it is easy to see that $v_n$ is still a maximizing sequence, since $\|e^{ith(D)}v_n\|_{L^r_{t,x}}=\|e^{ith(D)}w_n\|_{L^r_{t,x}}$. 

By Proposition \ref{prop:CC}, in order to conclude the proof it is sufficient to prove that \eqref{eq:4000} holds on $v_n$, up to subsequences, with $T=e^{ith(D)}$. Let us fix $t\in\R$; by the continuity of $e^{ith(D)}$ in $\dot H^s$, we have
\begin{equation*}
  e^{ith(D)}v_n\rightharpoonup e^{ith(D)}v
\end{equation*}
weakly in $\dot H^s$, as $n\to\infty$. Then, by the Rellich Theorem, for any $R>0$ 
\begin{equation}\label{eq:rellich}
  e^{ith(D)}v_n\to e^{ith(D)}v
\end{equation}
strongly in $L^2(B(R))$ where $B(R)=\{x\in \R^d:|x|< R\}$, as $n\to\infty$. Denote by
\begin{equation*}
  F_n(t):=\int_{|x|<R}\left|e^{ith(D)}(v_n-v)\right|^2dx=\left\|e^{ith(D)}(v_n
  -v)\right\|_{L^2(B(R))}^2.
\end{equation*}
By the H\"older inequality and Sobolev embedding we obtain
\begin{equation}\label{eq:equi}
  F_n(t)\leq CR^\frac{2s}{d}\left\|e^{ith(D)}(v_n-v)\right\|_{\dot H^s}^2\leq
  2CR^\frac{2s}{d};
\end{equation}
consequently, by \eqref{eq:rellich}, the Fubini and the Lebesgue Theorems we have that
\begin{equation*}
  \int_{-R}^RF_n(t)\,dt=\int_{-R}^R\int_{|x|<R}\left|e^{ith(D)}(v_n-v)\right|^2dx\,dt
  \to 0
\end{equation*}
as $n\to\infty$; this implies that, up to a subsequence, 
\begin{equation*}
  e^{ith(D)}(v_n-v)\to0 \quad\text{a.e. in } B(R)\times (-R, R).
\end{equation*}
The extraction of the subsequence depends on $R$; now repeat the argument on a discrete sequence of radii $R_n$ such that $R_n\to\infty$, as $n\to\infty$ and conclude, by a diagonal argument, that there exists a subsequence of $v_{n_k}$ of $v_n$ such that
\begin{equation}\label{eq:fine}
  e^{ith(D)}(v_n-v)\to0 \quad\text{a.e. in }\R\times \R^d.
\end{equation}
This, together with \eqref{eq:50} and Proposition \ref{prop:CC}, concludes the proof.

\end{document}